\documentclass[]{amsart}
\usepackage{amsmath, amsthm, amscd, amsfonts, amssymb, graphicx, color}
\usepackage[bookmarksnumbered, colorlinks, plainpages]{hyperref}

\theoremstyle{definition}

\theoremstyle{remark}

\numberwithin{equation}{section}

\begin{document}
\setcounter{page}{1}
\begin{center}
{\bf FACTORIZATION PROPERTIES AND TOPOLOGICAL CENTERS OF  LEFT MODULE ACTIONS }
\end{center}

\title[]{}

\author[]{KAZEM HAGHNEJAD AZAR   }

\address{}

\dedicatory{}



\begin{abstract}
 For a Banach left  module action, we will extend some propositions from Lau and $\ddot{U}$lger and others into general situations  and we establish  the relationships between  topological centers of the left module action with  the multiplier  and  factorization properties of left module actions.   We have some applications in  the dual groups.
 \end{abstract} \maketitle

\begin{center}
{\bf  1.Introduction and Preliminaries}
\end{center}

\noindent As is well-known [1], the second dual $A^{**}$ of $A$ endowed with the either Arens multiplications is a Banach algebra. The constructions of the two Arens multiplications in $A^{**}$ lead us to definition of topological centers for $A^{**}$ with respect both Arens multiplications. The topological centers of Banach algebras, module actions and applications of them  were introduced and discussed in [6, 8, 13, 14, 15, 16, 17, 21, 22], and they have attracted by some attentions.

\noindent Now we introduce some notations and definitions that we used
throughout  this paper.\\
 Let $A$ be a Banach algebra. We
say that a  net $(e_{\alpha})_{{\alpha}\in I}$ in $A$ is a left
approximate identity $(=LAI)$ [resp. right
approximate identity $(=RAI)$] if,
 for each $a\in A$,   $e_{\alpha}a\longrightarrow a$ [resp. $ae_{\alpha}\longrightarrow a$]. For $a\in A$
 and $a^\prime\in A^*$, we denote by $a^\prime a$
 and $a a^\prime$ respectively, the functionals on $A^*$ defined by $<a^\prime a,b>=<a^\prime,ab>=a^\prime(ab)$ and $<a a^\prime,b>=<a^\prime,ba>=a^\prime(ba)$ for all $b\in A$.
  The Banach algebra $A$ is embedded in its second dual via the identification
 $<a,a^\prime>$ - $<a^\prime,a>$ for every $a\in
A$ and $a^\prime\in
A^*$. We denote the set   $\{a^\prime a:~a\in A~ and ~a^\prime\in
  A^*\}$ and
  $\{a a^\prime:~a\in A ~and ~a^\prime\in A^*\}$ by $A^*A$ and $AA^*$, respectively, clearly these two sets are subsets of $A^*$.  Let $A$ has a $BAI$. If the
equality $A^*A=A^*,~~(AA^*=A^*)$ holds, then we say that $A^*$
factors on the left (right). If both equalities $A^*A=AA^*=A^*$
hold, then we say
that $A^*$  factors on both sides.
 Let $X,Y,Z$ be normed spaces and $m:X\times Y\rightarrow Z$ be a bounded bilinear mapping. Arens in [1] offers two natural extensions $m^{***}$ and $m^{t***t}$ of $m$ from $X^{**}\times Y^{**}$ into $Z^{**}$ as following:\\
 \noindent1. $m^*:Z^*\times X\rightarrow Y^*$,~~~~~given by~~~$<m^*(z^\prime,x),y>=<z^\prime, m(x,y)>$ ~where $x\in X$, $y\in Y$, $z^\prime\in Z^*$,\\
 2. $m^{**}:Y^{**}\times Z^{*}\rightarrow X^*$,~~given by $<m^{**}(y^{\prime\prime},z^\prime),x>=<y^{\prime\prime},m^*(z^\prime,x)>$ ~where $x\in X$, $y^{\prime\prime}\in Y^{**}$, $z^\prime\in Z^*$,\\
3. $m^{***}:X^{**}\times Y^{**}\rightarrow Z^{**}$,~ given by~ ~ ~$<m^{***}(x^{\prime\prime},y^{\prime\prime}),z^\prime>$ $=<x^{\prime\prime},m^{**}(y^{\prime\prime},z^\prime)>$ \\~where ~$x^{\prime\prime}\in X^{**}$, $y^{\prime\prime}\in Y^{**}$, $z^\prime\in Z^*$.\\
The mapping $m^{***}$ is the unique extension of $m$ such that $x^{\prime\prime}\rightarrow m^{***}(x^{\prime\prime},y^{\prime\prime})$ from $X^{**}$ into $Z^{**}$ is $weak^*-to-weak^*$ continuous for every $y^{\prime\prime}\in Y^{**}$, but the mapping $y^{\prime\prime}\rightarrow m^{***}(x^{\prime\prime},y^{\prime\prime})$ is not in general $weak^*-to-weak^*$ continuous from $Y^{**}$ into $Z^{**}$ unless $x^{\prime\prime}\in X$. Hence the first topological center of $m$ may  be defined as following
$$Z_1(m)=\{x^{\prime\prime}\in X^{**}:~~y^{\prime\prime}\rightarrow m^{***}(x^{\prime\prime},y^{\prime\prime})~~is~~weak^*-to-weak^*-continuous\}.$$
Let now $m^t:Y\times X\rightarrow Z$ be the transpose of $m$ defined by $m^t(y,x)=m(x,y)$ for every $x\in X$ and $y\in Y$. Then $m^t$ is a continuous bilinear map from $Y\times X$ to $Z$, and so it may be extended as above to $m^{t***}:Y^{**}\times X^{**}\rightarrow Z^{**}$.
 The mapping $m^{t***t}:X^{**}\times Y^{**}\rightarrow Z^{**}$ in general is not equal to $m^{***}$, see [1], if $m^{***}=m^{t***t}$, then $m$ is called Arens regular. The mapping $y^{\prime\prime}\rightarrow m^{t***t}(x^{\prime\prime},y^{\prime\prime})$ is $weak^*-to-weak^*$ continuous for every $y^{\prime\prime}\in Y^{**}$, but the mapping $x^{\prime\prime}\rightarrow m^{t***t}(x^{\prime\prime},y^{\prime\prime})$ from $X^{**}$ into $Z^{**}$ is not in general  $weak^*-to-weak^*$ continuous for every $y^{\prime\prime}\in Y^{**}$. So we define the second topological center of $m$ as
$$Z_2(m)=\{y^{\prime\prime}\in Y^{**}:~~x^{\prime\prime}\rightarrow m^{t***t}(x^{\prime\prime},y^{\prime\prime})~~is~~weak^*-to-weak^*-continuous\}.$$
It is clear that $m$ is Arens regular if and only if $Z_1(m)=X^{**}$ or $Z_2(m)=Y^{**}$. Arens regularity of $m$ is equivalent to the following
$$\lim_i\lim_j<z^\prime,m(x_i,y_j)>=\lim_j\lim_i<z^\prime,m(x_i,y_j)>,$$
whenever both limits exist for all bounded sequences $(x_i)_i\subseteq X$ , $(y_i)_i\subseteq Y$ and $z^\prime\in Z^*$, see [6, 18].\\
 The regularity of a normed algebra $A$ is defined to be the regularity of its algebra multiplication when considered as a bilinear mapping. Let $a^{\prime\prime}$ and $b^{\prime\prime}$ be elements of $A^{**}$, the second dual of $A$. By $Goldstin^,s$ Theorem [6, P.424-425], there are nets $(a_{\alpha})_{\alpha}$ and $(b_{\beta})_{\beta}$ in $A$ such that $a^{\prime\prime}=weak^*-\lim_{\alpha}a_{\alpha}$ ~and~  $b^{\prime\prime}=weak^*-\lim_{\beta}b_{\beta}$. So it is easy to see that for all $a^\prime\in A^*$,
$$\lim_{\alpha}\lim_{\beta}<a^\prime,m(a_{\alpha},b_{\beta})>=<a^{\prime\prime}b^{\prime\prime},a^\prime>$$ and
$$\lim_{\beta}\lim_{\alpha}<a^\prime,m(a_{\alpha},b_{\beta})>=<a^{\prime\prime}ob^{\prime\prime},a^\prime>,$$
where $a^{\prime\prime}b^{\prime\prime}$ and $a^{\prime\prime}ob^{\prime\prime}$ are the first and second Arens products of $A^{**}$, respectively, see [6, 14, 18].\\
The mapping $m$ is left strongly Arens irregular if $Z_1(m)=X$ and $m$ is right strongly Arens irregular if $Z_2(m)=Y$.\\
Let $A$ and $B$ be normed spaces. In this paper, if $T$ is a continuous linear operator from $A$ into $B$, then we write $T\in \mathbf{B}(A,B)$.\\
This paper is organized as follows.\\
{\bf A.} Let $B$ be a Banach left $A-module$ and $(e_{\alpha})_{\alpha}\subseteq A$ be a $LBAI$ for $B$. Then the following assertions hold.
\begin{enumerate}
\item For each $b^\prime\in B^*$, $\pi_\ell^*(b^\prime,e_\alpha)\stackrel{w^*} {\rightarrow}b^\prime$.
\item $B^*$ factors on the left with respect to $A$ if and only if $B^{**}$ has a $W^*LBAI$ $(e_{\alpha})_{\alpha}\subseteq A$.
\item $B^{**}$ has a $W^*LBAI$ $(e_{\alpha})_{\alpha}\subseteq A$ if and only if $B^{**}$ has a left unit element $e^{\prime \prime}\in A^{**}$ such that $e_{\alpha}\stackrel{w^*} {\rightarrow}e^{\prime \prime}$.
\end{enumerate}
{\bf B.} Let $B$ be a Banach left $A-module$ and suppose that $b^\prime\in wap_\ell(B)$. Let  $a^{\prime\prime}\in A^{**}$ and
$(a_{\alpha})_{\alpha}\subseteq A$ such that  $a_{\alpha} \stackrel{w^*} {\rightarrow}a^{\prime\prime}$ in $A^{**}$. Then we have
$$\pi^*_\ell(b^\prime,a_\alpha) \stackrel{w} {\rightarrow}\pi_\ell^{****}(b^\prime,a^{\prime\prime}).$$
{\bf C.} Let $B$ be a Banach left $A-module$ and $B^*$ factors on the left with respect to $A$. If $AA^{**}\subseteq  Z_{B^{**}}(A^{**})$, then  $Z_{B^{**}}(A^{**})=A^{**}$.\\
{\bf D.} Let $B$ be a Banach left $A-module$ and $B^{**}$ has a $LBAI$ with respect to $A^{**}$. Then $B^{**}$ has a left unit with respect to $A^{**}$.\\
{\bf E.} Let $B$ be a Banach left $A-module$ and it has a $LBAI$ with respect to $A$. Then we have the following assertions.
\begin{enumerate}
\item  $B^*$ factors on the left with respect to $A$ if and only if for each $b^\prime\in B^*$, we have $\pi^*_\ell(b^\prime,e_\alpha) \stackrel{w} {\rightarrow}b^\prime$ in $B^*$.
\item $B$ factors on the left with respect to $A$ if and only if for each $b\in B$, we have $\pi^*_\ell(b,e_\alpha) \stackrel{w} {\rightarrow}b$ in $B$.
\end{enumerate}
{\bf F.} Let $B$ be a Banach left $A-module$ and $A$ has a $LBAI$ $(e_{\alpha})_{\alpha}\subseteq A$ such that $e_{\alpha} \stackrel{w^*} {\rightarrow}e^{\prime\prime}$ in $A^{**}$ where $e^{\prime\prime}$ is a left unit for $A^{**}$. Suppose that $Z^t_{e^{\prime\prime}}(B^{**})=B^{**}$. {Then, $B$ factors on the right with respect to $A$ if and only if
$e^{\prime\prime}$ is a left unit for $B^{**}$.\\
{\bf G.} Let $B$ be a left  Banach $A-module$ and $T\in \mathbf{B}(A,B)$. Consider the following statements.
\begin{enumerate}
\item $T=\ell_b$, for some $b\in B$.
\item  $T^{**}(a^{\prime\prime})=\pi_\ell^{***}(a^{\prime\prime},b^{\prime\prime})$ for some $b^{\prime\prime}\in B^{**}$ such that $\widetilde{Z}_{b^{\prime\prime}}(A^{**})=A^{**}$.
\item $T^*(B^*)\subseteq BB^*$.
\end{enumerate}
Then $(1)\Rightarrow (2)\Rightarrow (3)$.\\
If we take $T\in RM(A,B)$ and if $B$ has a sequential $BAI$ and it is $WSC$, then (1), (2) and (3) are equivalent.\\\\

\begin{center}
\textbf{{ 2. Factorization  properties and  topological centers of  left module actions}}
\end{center}

\noindent Let $B$ be a Banach $A-bimodule$, and let\\
$$\pi_\ell:~A\times B\rightarrow B~~~and~~~\pi_r:~B\times A\rightarrow B.$$
be the left and right module actions of $A$ on $B$. Then $B^{**}$ is a Banach $A^{**}-bimodule$ with module actions
$$\pi_\ell^{***}:~A^{**}\times B^{**}\rightarrow B^{**}~~~and~~~\pi_r^{***}:~B^{**}\times A^{**}\rightarrow B^{**}.$$
Similarly, $B^{**}$ is a Banach $A^{**}-bimodule$ with module actions\\
$$\pi_\ell^{t***t}:~A^{**}\times B^{**}\rightarrow B^{**}~~~and~~~\pi_r^{t***t}:~B^{**}\times A^{**}\rightarrow B^{**}.$$
We may therefore define the topological centers of the right and left module actions of $A$ on $B$ as follows:\\
$$Z_{A^{**}}(B^{**})=Z(\pi_r)=\{b^{\prime\prime}\in B^{**}:~the~map~~a^{\prime\prime}\rightarrow \pi_r^{***}(b^{\prime\prime}, a^{\prime\prime})~:~A^{**}\rightarrow B^{**}$$$$~is~~~weak^*-to-weak^*~continuous\}$$
$$Z_{B^{**}}(A^{**})=Z(\pi_\ell)=\{a^{\prime\prime}\in A^{**}:~the~map~~b^{\prime\prime}\rightarrow \pi_\ell^{***}(a^{\prime\prime}, b^{\prime\prime})~:~B^{**}\rightarrow B^{**}$$$$~is~~~weak^*-to-weak^*~continuous\}$$
$$Z_{A^{**}}^t(B^{**})=Z(\pi_\ell^t)=\{b^{\prime\prime}\in B^{**}:~the~map~~a^{\prime\prime}\rightarrow \pi_\ell^{t***}(b^{\prime\prime}, a^{\prime\prime})~:~A^{**}\rightarrow B^{**}$$$$~is~~~weak^*-to-weak^*~continuous\}$$
$$Z_{B^{**}}^t(A^{**})=Z(\pi_r^t)=\{a^{\prime\prime}\in A^{**}:~the~map~~b^{\prime\prime}\rightarrow \pi_r^{t***}(a^{\prime\prime}, b^{\prime\prime})~:~B^{**}\rightarrow B^{**}$$$$~is~~~weak^*-to-weak^*~continuous\}$$
We note also that if $B$ is a left (resp. right) Banach $A-module$ and $\pi_\ell:~A\times B\rightarrow B$~(resp. $\pi_r:~B\times A\rightarrow B$) is left (resp. right) module action of $A$ on $B$, then $B^*$ is a right (resp. left) Banach $A-module$. \\
We write $ab=\pi_\ell(a,b)$, $ba=\pi_r(b,a)$, $\pi_\ell(a_1a_2,b)=\pi_\ell(a_1,a_2b)$,  $\pi_r(b,a_1a_2)=\pi_r(ba_1,a_2)$,~\\
$\pi_\ell^*(a_1b^\prime, a_2)=\pi_\ell^*(b^\prime, a_2a_1)$,~
$\pi_r^*(b^\prime a, b)=\pi_r^*(b^\prime, ab)$,~ for all $a_1,a_2, a\in A$, $b\in B$ and  $b^\prime\in B^*$
when there is no confusion.\\
Let $B$ be a left Banach $A-module$ and  $e$ be a left  unit element of $A$. Then we say that $e$ is a left unit (resp. weakly left unit)  $A-module$ for $B$, if $\pi_\ell(e,b)=b$ (resp. $<b^\prime , \pi_\ell(e,b)>=<b^\prime , b>$ for all $b^\prime\in B^*$) where $b\in B$. The definition of right unit (resp. weakly right unit) $A-module$ is similar.\\
We say that a Banach $A-bimodule$ $B$ is a unital $A-module$, if $B$ has left and right unit $A-module$ that are equal then we say that $B$ is unital $A-module$.\\
Let $B$ be a left Banach $A-module$ and $(e_{\alpha})_{\alpha}\subseteq A$ be a LAI [resp. weakly left approximate identity(=WLAI)] for $A$. We say that $(e_{\alpha})_{\alpha}$ is left approximate identity  ($=LAI$)[ resp. weakly left approximate identity  (=$WLAI$)] for $B$, if for all $b\in B$, we have $\pi_\ell (e_{\alpha},b) \stackrel{} {\rightarrow}
b$ ( resp. $\pi_\ell (e_{\alpha},b) \stackrel{w} {\rightarrow}
b$). The definition of the right approximate identity ($=RAI$)[ resp. weakly right approximate identity ($=WRAI$)] is similar.\\
We say that $(e_{\alpha})_{\alpha}$ is a approximate identity  ($=AI$)[ resp. weakly approximate identity  ($WAI$)] for $B$, if $B$ has the same left and right approximate identity  [ resp. weakly left and right approximate identity ].\\
Let $(e_{\alpha})_{\alpha}\subseteq A$ be $weak^*$ left approximate identity for $A^{**}$. We say that $(e_{\alpha})_{\alpha}$ is $weak^*$ left approximate identity $A^{**}-module$ ($=W^*LAI~A^{**}-module$) for $B^{**}$, if for all $b^{\prime\prime}\in B^{**}$, we have $\pi_\ell^{***} (e_{\alpha},b^{\prime\prime}) \stackrel{w^*} {\rightarrow}
b^{\prime\prime}$. The definition of the $weak^*$ right approximate identity $A^{**}-module$ ($=W^*RAI~A^{**}-module$) is similar.\\
 We say that $(e_{\alpha})_{\alpha}$ is a $weak^*$ approximate identity $A^{**}-module$ ($=W^*AI~A^{**}-module$) for $B^{**}$, if $B^{**}$ has the same $weak^*$ left and right approximate identity $A^{**}-module$.\\
Let $B$ be a Banach $A-bimodule$. We say that $B$ is a left [resp. right] factors with respect to $A$, if $BA=B$ [resp. $AB=B$].\\\\

\noindent{\it{\bf Theorem 2-1.}} Let $B$ be a Banach left $A-module$ and $(e_{\alpha})_{\alpha}\subseteq A$ be a $LBAI$ for $B$. Then the following assertions hold.
\begin{enumerate}
\item For each $b^\prime\in B^*$, we have  $\pi_\ell^*(b^\prime,e_\alpha)\stackrel{w^*} {\rightarrow}b^\prime$.
\item $B^*$ factors on the left with respect to $A$ if and only if $B^{**}$ has a $W^*LBAI$ $(e_{\alpha})_{\alpha}\subseteq A$.
\item $B^{**}$ has a $W^*LBAI$ $(e_{\alpha})_{\alpha}\subseteq A$ if and only if $B^{**}$ has a left unit element $e^{\prime \prime}\in A^{**}$ such that $e_{\alpha}\stackrel{w^*} {\rightarrow}e^{\prime \prime}$.
\end{enumerate}
\begin{proof}
\begin{enumerate}
\item Let $b\in B$ and $b^\prime\in B^*$. Since $\mid <b^\prime , \pi_\ell(e_\alpha ,b)\mid\leq \parallel b^\prime\parallel\parallel \pi_\ell(e_\alpha ,b)\parallel$, we have the following equality
$$\lim_\alpha <\pi^*_\ell(b^\prime,e_\alpha),b>=\lim_\alpha<b^\prime,\pi_\ell(e_\alpha,b)>= 0.$$
It follows that  $\pi^*_\ell(b^\prime,e_\alpha)\stackrel{w^*} {\rightarrow}0$.

\item  Let $B^*$ factors on the left with respect to $A$. Then for every $b^\prime\in B^*$, there are $x^\prime\in B^*$ and $a\in A$ such that $b^\prime=x^\prime a$. Then for every $b^{\prime\prime}\in B^{**}$, we have
$$<\pi^{***}_\ell(e_\alpha,b^{\prime\prime}),b^\prime >= <e_\alpha,\pi^{**}_\ell(b^{\prime\prime},b^\prime )>=
<\pi^{**}_\ell(b^{\prime\prime},b^\prime ),e_\alpha>$$$$=<b^{\prime\prime},\pi^{*}_\ell(b^\prime ,e_\alpha)>=
<b^{\prime\prime},\pi^{*}_\ell(x^\prime a,e_\alpha)>= <b^{\prime\prime},\pi^{*}_\ell(x^\prime ,ae_\alpha)>$$$$=
 <\pi^{**}_\ell(b^{\prime\prime},x^\prime ),ae_\alpha>\rightarrow   <\pi^{**}_\ell(b^{\prime\prime},x^\prime ),a>
 $$$$=< b^{\prime\prime},b^\prime >.$$
It follows that
$$\pi^{***}_\ell(e_\alpha,b^{\prime\prime})\stackrel{w^*} {\rightarrow}b^{\prime\prime},$$
 and so $B^{**}$ has $W^*LBAI$.\\
Conversely, let $b^\prime\in B^*$. Then for every $b^{\prime\prime}\in B^{**}$, we have
$$<b^{\prime\prime},\pi^{*}_\ell(b^\prime ,e_\alpha)>=<\pi^{***}_\ell(e_\alpha,b^{\prime\prime}),b^\prime >\rightarrow< b^{\prime\prime},b^\prime >.$$
It follows that
$$\pi^{*}_\ell(b^\prime ,e_\alpha)\stackrel{w} {\rightarrow} b^\prime ,$$
and so by Cohen factorization theorem, we are done.

\item Assume that  $B^{**}$ has a $W^*LBAI$ $(e_{\alpha})_{\alpha}\subseteq A$. Without loss generality, let $e^{\prime\prime}\in A^{**}$ be a left unit for $A^{**}$ with respect to the first Arens product such that $e_{\alpha}\stackrel{w^*} {\rightarrow}e^{\prime\prime}$. Then for each $b^\prime\in B^*$, we have

 $$<\pi^{***}_\ell(e^{\prime\prime},b^{\prime\prime}),b^\prime >=<e^{\prime\prime},\pi^{**}_\ell(b^{\prime\prime},b^\prime ) >$$$$=\lim_\alpha<e_{\alpha},\pi^{**}_\ell(b^{\prime\prime},b^\prime) >=\lim_\alpha<\pi^{**}_\ell(b^{\prime\prime},b^\prime ),e_{\alpha}>
 $$$$=\lim_\alpha<b^{\prime\prime},\pi^{*}_\ell(b^\prime ,e_{\alpha})>
=\lim_\alpha<\pi^{****}_\ell(b^\prime ,e_{\alpha}),b^{\prime\prime}>
$$$$=\lim_\alpha<b^\prime ,\pi^{***}_\ell(e_{\alpha},b^{\prime\prime})>
=\lim_\alpha<\pi^{***}_\ell(e_{\alpha},b^{\prime\prime}),b^\prime >$$$$= < b^{\prime\prime},b^\prime >.$$
Thus $e^{\prime\prime}\in A^{**}$ is a left unit for $B^{**}$.\\
Conversely, let  $e^{\prime\prime}\in A^{**}$ be a left unit for $B^{**}$ and assume that   $e_{\alpha}\stackrel{w^*} {\rightarrow}e^{\prime\prime}$ in $A^{**}$. Then for  every $b^{\prime\prime}\in B^{**}$ and $b^\prime\in B^*$, we have
$$<\pi^{***}_\ell(e_\alpha,b^{\prime\prime}),b^\prime >=<e_\alpha,\pi^{**}_\ell(b^{\prime\prime},b^\prime )>$$$$\rightarrow
<e^{\prime\prime},\pi^{**}_\ell(b^{\prime\prime},b^\prime )>=<\pi^{***}_\ell(e^{\prime\prime},b^{\prime\prime}),b^\prime >$$
$$= < b^{\prime\prime},b^\prime >.$$
It follows that $\pi^{***}_\ell(e_\alpha,b^{\prime\prime})\stackrel{w^*} {\rightarrow}b^{\prime\prime}.$\\
\end{enumerate}
\end{proof}

\noindent{\it{\bf Corollary 2-2.}} Let $B$ be a Banach left $A-module$ and $A$ has a $BLAI$. If $B^{**}$ has a $W^*LBAI$ , then
$$\{a^{\prime\prime}\in A^{**}:~Aa^{\prime\prime}\subseteq A\}\subseteq Z_{B^{**}}(A^{**}).$$
\begin{proof} By using the proceeding theorem, since $B^{**}$ has $W^*LBAI$, $B^*$ factors on the left with respect to $A$. Suppose that $b^\prime\in B^*$. Then there are $x^\prime\in B^*$ and $a\in A$ such that $b^\prime=x^\prime a$. Assume that $a^{\prime\prime}\in A^{**}$ such that $Aa^{\prime\prime}\subseteq A$. Let $b^{\prime\prime}\in B^{**}$ and $(b^{\prime\prime}_\alpha)_\alpha\subseteq B^{**}$ such that $b^{\prime\prime}_\alpha\stackrel{w^*} {\rightarrow}b^{\prime\prime}$ in $B^{**}$. Then we have the following equality
$$ \lim_\alpha<\pi_\ell^{***}(a^{\prime\prime},b^{\prime\prime}_\alpha),b^\prime>
=\lim_\alpha<\pi_\ell^{***}(a^{\prime\prime},b^{\prime\prime}_\alpha),x^\prime a >$$$$=\lim_\alpha<a\pi_\ell^{***}(a^{\prime\prime},b^{\prime\prime}_\alpha),x^\prime >=
<\pi_\ell^{***}(aa^{\prime\prime},b^{\prime\prime}),x^\prime >=
<\pi_\ell^{***}(a^{\prime\prime},b^{\prime\prime}),b^\prime >.$$

It follows that $a^{\prime\prime}\in Z_{B^{**}}(A^{**}).$\\
\end{proof}

In the proceeding corollary, if we take $B=A$, then we have the following conclusion
$$\{a^{\prime\prime}\in A^{**}:~Aa^{\prime\prime}\subseteq A\}\subseteq Z_{1}(A^{**}).$$\\

A functional $a^\prime$ in $A^*$ is said to be $wap$ (weakly almost
 periodic) on $A$ if the mapping $a\rightarrow a^\prime a$ from $A$ into
 $A^{*}$ is weakly compact. The proceeding definition is equivalent to the following condition, see [6, 14, 18].\\
 For any two net $(a_{\alpha})_{\alpha}$ and $(b_{\beta})_{\beta}$
 in $\{a\in A:~\parallel a\parallel\leq 1\}$, we have\\
$$\\lim_{\alpha}\\lim_{\beta}<a^\prime,a_{\alpha}b_{\beta}>=\\lim_{\beta}\\lim_{\alpha}<a^\prime,a_{\alpha}b_{\beta}>,$$
whenever both iterated limits exist. The collection of all $wap$
functionals on $A$ is denoted by $wap(A)$. Also we have
$a^{\prime}\in wap(A)$ if and only if $<a^{\prime\prime}b^{\prime\prime},a^\prime>=<a^{\prime\prime}ob^{\prime\prime},a^\prime>$ for every $a^{\prime\prime},~b^{\prime\prime} \in
A^{**}$.\\
Let $B$ be a Banach left $A-module$. Then, $b^\prime\in B^*$ is said to be left weakly almost periodic functional if the set $\{\pi_\ell(b^\prime,a):~a\in A,~\parallel a\parallel\leq 1\}$ is relatively weakly compact. We denote by $wap_\ell(B)$ the closed subspace of $B^*$ consisting of all the left weakly almost periodic functionals in $B^*$.\\
The definition of the right weakly almost periodic functional ($=wap_r(B)$) is the same.\\
By  [18],  $b^\prime\in wap_\ell(B)$ if and only if $$<\pi_\ell^{***}(a^{\prime\prime},b^{\prime\prime}),b^\prime>=
<\pi_\ell^{t***t}(a^{\prime\prime},b^{\prime\prime}),b^\prime>$$
for all $a^{\prime\prime}\in A^{**}$ and $b^{\prime\prime}\in B^{**}$.
Thus, we can write \\
$$wap_\ell(B)=\{ b^\prime\in B^*:~<\pi_\ell^{***}(a^{\prime\prime},b^{\prime\prime}),b^\prime>=
<\pi_\ell^{t***t}(a^{\prime\prime},b^{\prime\prime}),b^\prime>~~$$$$for~~all~~a^{\prime\prime}\in A^{**},~b^{\prime\prime}\in B^{**}\}.$$\\

\noindent{\it{\bf Theorem 2-3.}} Let $B$ be a Banach left $A-module$ and suppose that $b^\prime\in wap_\ell(B)$. Let  $a^{\prime\prime}\in A^{**}$ and
$(a_{\alpha})_{\alpha}\subseteq A$ such that  $a_{\alpha} \stackrel{w^*} {\rightarrow}a^{\prime\prime}$ in $A^{**}$. Then we have
$$\pi^*_\ell(b^\prime,a_\alpha) \stackrel{w} {\rightarrow}\pi_\ell^{****}(b^\prime,a^{\prime\prime}).$$
\begin{proof}
Assume that $b^{\prime\prime}\in B^{**}$. Then we have the following equality

$$<\pi_\ell^{****}(b^\prime,a^{\prime\prime}),b^{\prime\prime}>=
<\pi_\ell^{***}(a^{\prime\prime},b^{\prime\prime}),b^\prime>=
\lim_\alpha<\pi_\ell^{***}(a_{\alpha},b^{\prime\prime}),b^\prime>
$$$$=\lim_\alpha<b^{\prime\prime},\pi_\ell^{*}(b^\prime, a_{\alpha})>.$$
Now suppose that $(b_{\beta}^{\prime\prime})_{\beta}\subseteq B^{**}$ such that  $b^{\prime\prime}_{\beta} \stackrel{w^*} {\rightarrow}b^{\prime\prime}$. Since $b^\prime\in wap_\ell(B)$, we have
$$<\pi_\ell^{****}(b^\prime,a^{\prime\prime}),b_\beta^{\prime\prime}>=
<\pi_\ell^{***}(a^{\prime\prime},b_\beta^{\prime\prime}),b^\prime>\rightarrow
<\pi_\ell^{***}(a^{\prime\prime},b^{\prime\prime}),b^\prime>
$$$$=<\pi_\ell^{****}(b^\prime,a^{\prime\prime}),b^{\prime\prime}>.$$
Thus $\pi_\ell^{****}(b^\prime,a^{\prime\prime})\in (B^{**},weak^*)^*=B^*$.
So we conclude that
$$\pi^*_\ell(b^\prime,a_\alpha) \stackrel{w} {\rightarrow}\pi_\ell^{****}(b^\prime,a^{\prime\prime})~in ~ B^{**}.$$\\
\end{proof}

In the proceeding corollary, if we take $B=A$, then we obtain the following result.\\
Suppose that $a^\prime\in wap(A)$ and $a^{\prime\prime}\in A^{**}$ such that
$a_{\alpha} \stackrel{w^*} {\rightarrow} a^{\prime\prime}$ where\\
$(a_{\alpha})_{\alpha}\subseteq A$. Then we have $a^\prime a_{\alpha} \stackrel{w}
{\rightarrow} a^\prime a^{\prime\prime}$.\\\\

\noindent{\it{\bf Theorem 2-4.}} Let $B$ be a Banach left $A-module$ and it has a $BLAI$ $(e_{\alpha})_{\alpha}\subseteq A$. Suppose that $b^\prime\in wap_\ell(B)$. Then we have
$$\pi^*_\ell(b^\prime,e_\alpha) \stackrel{w} {\rightarrow}b^\prime.$$
\begin{proof}
Let $b^{\prime\prime}\in B^{**}$ and $(b_\beta)_\beta\subseteq B$ such that $b_\beta\stackrel{w^*} {\rightarrow}b^{\prime\prime}$ in $B^{**}$. Then for every $b^\prime\in wap_\ell(B)$, we have the following equality
$$\lim_\alpha<b^{\prime\prime}, \pi^*_\ell(b^\prime,e_\alpha)>=
\lim_\alpha<\pi^{****}_\ell( b^\prime,e_\alpha),b^{\prime\prime}>$$$$=
\lim_\alpha< b^\prime,\pi^{***}_\ell(e_\alpha,b^{\prime\prime})>=
\lim_\alpha< \pi^{***}_\ell(e_\alpha,b^{\prime\prime}),b^\prime>$$$$=
\lim_\alpha\lim_\beta< \pi_\ell(e_\alpha,b_\beta),b^\prime>=
\lim_\beta\lim_\alpha< \pi_\ell(e_\alpha,b_\beta),b^\prime>$$$$=
\lim_\beta< b_\beta,b^\prime>=< b^{\prime\prime},b^\prime>.$$
It follows that
$$\pi^*_\ell(b^\prime,e_\alpha) \stackrel{w} {\rightarrow}b^\prime.$$
\end{proof}

\noindent{\it{\bf Corollary 2-5.}} Let $B$ be a Banach left $A-module$ and it has a $BLAI$ $(e_{\alpha})_{\alpha}\subseteq A$. Suppose that $ wap_\ell(B)=B^*$. Then $B^*$ factors on the left with respect to $A$.\\\\

\noindent{\it{\bf Corollary 2-6.}} Let $A$ be an Arens regular Banach algebra with $LBAI$. Then $A^*$ factors on the left.\\\\

\noindent{\it{\bf Example 2-7.}} i) Let $G$ be finite group. Then we have the following equality
$$M(G)^*L^1(G)=M(G)^*~and ~L^\infty (G)L^1(G)=L^\infty (G).$$
ii) Consider the Banach algebra $(\ell^1,.)$ that is Arens regular Banach algebra with unit element. Then we have $\ell^\infty .\ell^1=\ell^\infty$.\\\\

\noindent{\it{\bf Theorem 2-8.}} Let $B$ be a Banach left $A-module$ and $B^*$ factors on the left with respect to $A$. If $AA^{**}\subseteq  Z_{B^{**}}(A^{**})$, then  $Z_{B^{**}}(A^{**})=A^{**}$.
\begin{proof}
Let $b^{\prime\prime}\in B^{**}$ and $(b^{\prime\prime}_\alpha)_\alpha\subseteq B^{**}$ such that $b^{\prime\prime}_\alpha\stackrel{w^*} {\rightarrow}b^{\prime\prime}$ in $B^{**}$. Suppose that $a^{\prime\prime}\in A^{**}$. Since $B^*$ factors on the left with respect to $A$,  for every $b^\prime\in B^*$, there are $x^\prime\in B^*$ and $a\in A$ such that $b^\prime=x^\prime a$. Since $aa^{\prime\prime}\in Z_{B^{**}}(A^{**})$, we have
$$<\pi^{***}_\ell(a^{\prime\prime},b_\alpha^{\prime\prime}),b^\prime >=
<\pi^{***}_\ell(a^{\prime\prime},b_\alpha^{\prime\prime}),x^\prime a>$$$$=
<a\pi^{***}_\ell(a^{\prime\prime},b_\alpha^{\prime\prime}),x^\prime >=
<\pi^{***}_\ell(aa^{\prime\prime},b_\alpha^{\prime\prime}),x^\prime >$$$$
\rightarrow <\pi^{***}_\ell(aa^{\prime\prime},b^{\prime\prime}),x^\prime >
=<\pi^{***}_\ell(a^{\prime\prime},b^{\prime\prime}),b^\prime >.$$
It follows that
$$\pi^{***}_\ell(a^{\prime\prime},b_\alpha^{\prime\prime})\stackrel{w^*} {\rightarrow}
<\pi^{***}_\ell(a^{\prime\prime},b^{\prime\prime}),$$
and so $a^{\prime\prime}\in Z_{B^{**}}(A^{**})$.\\
\end{proof}

\noindent{\it{\bf Corollary 2-9.}} Let $A$ be a Banach algebra and $A^*$ factors on the left. If $AA^{**}\subseteq Z_1(A^{**})$, then $A$ is Arens regular.\\\\

\noindent{\it{\bf Theorem 2-10.}} Let $B$ be a Banach left $A-module$ and $B^{**}$ has a $LBAI$ with respect to $A^{**}$. Then $B^{**}$ has a left unit with respect to $A^{**}$.
\begin{proof}
Assume that $(e^{\prime\prime}_{\alpha})_{\alpha}\subseteq A^{**}$ is a $LBAI$ for $B^{**}$. By passing to a subnet, we may suppose that there is  $e^{\prime\prime}\in A^{**}$ such that $e^{\prime\prime}_{\alpha} \stackrel{w^*} {\rightarrow}e^{\prime\prime}$ in $A^{**}$. Then for every $b^{\prime\prime}\in B^{**}$ and $b^\prime\in B^*$, we have

$$<\pi_\ell^{***}(e^{\prime\prime},b^{\prime\prime}),b^\prime>=
<e^{\prime\prime},\pi_\ell^{**}(b^{\prime\prime},b^\prime)>=
\lim_\alpha <e_\alpha^{\prime\prime},\pi_\ell^{**}(b^{\prime\prime},b^\prime)>$$$$=
\lim_\alpha <\pi_\ell^{***}(e_\alpha^{\prime\prime},b^{\prime\prime}),b^\prime>=
<b^{\prime\prime},b^\prime>.$$
It follows that $\pi_\ell^{***}(e^{\prime\prime},b^{\prime\prime})=b^{\prime\prime}$.\\\end{proof}

\noindent{\it{\bf Corollary 2-11.}} Let $A$ be a Banach algebra and $A^{**}$ has a $LBAI$. Then $A^{**}$ has a left unit with respect to the first Arens product.\\\\

\noindent{\it{\bf Theorem 2-12.}} Let $B$ be a Banach left $A-module$ and it has a $LBAI$ with respect to $A$. Then we have the following assertions.
\begin{enumerate}
\item  $B^*$ factors on the left with respect to $A$ if and only if for each $b^\prime\in B^*$, we have $\pi^*_\ell(b^\prime,e_\alpha) \stackrel{w} {\rightarrow}b^\prime$ in $B^*$.
\item $B$ factors on the left with respect to $A$ if and only if for each $b\in B$, we have $\pi^*_\ell(b,e_\alpha) \stackrel{w} {\rightarrow}b$ in $B$.
\end{enumerate}
\begin{proof}\begin{enumerate}
\item Assume that  $B^*$ factors on the left with respect to $A$. Then  for every $b^\prime\in B^*$, there are $x^\prime\in B^*$ and $a\in A$ such that $b^\prime=x^\prime a$. Then for every $b^{\prime\prime}\in B^{**}$, we have
$$<b^{\prime\prime}, \pi^*_\ell(b^\prime,e_\alpha)>= <b^{\prime\prime}, \pi^*_\ell(x^\prime a,e_\alpha)>=<b^{\prime\prime}, \pi^*_\ell(x^\prime ,ae_\alpha)>$$$$=
<\pi^{**}_\ell(b^{\prime\prime}, x^\prime), ae_\alpha)>\rightarrow <\pi^{**}_\ell(b^{\prime\prime}, x^\prime), a>$$$$=
<b^{\prime\prime},b^\prime>.$$
It follows that
 $\pi^*_\ell(b^\prime,e_\alpha) \stackrel{w} {\rightarrow}b^\prime$ in $B^*$.
The converse by Cohen factorization theorem hold.

\item It is similar to the proceeding proof.\\
\end{enumerate}\end{proof}

In the proceeding theorem, if we take $B=A$, we obtain Lemma 2.1 from [14].\\

Let $B$ be a Banach  $A-bimodule$ and $a^{\prime\prime}\in A^{**}$. We define the locally topological centers of the left and right module actions of $a^{\prime\prime}$ on $B$, respectively, as follows\\
$$Z_{a^{\prime\prime}}^t(B^{**})=Z_{a^{\prime\prime}}^t(\pi_\ell^t)=\{b^{\prime\prime}\in B^{**}:~~~\pi^{t***t}_\ell(a^{\prime\prime},b^{\prime\prime})=
\pi^{***}_\ell(a^{\prime\prime},b^{\prime\prime})\},$$
$$Z_{a^{\prime\prime}}(B^{**})=Z_{a^{\prime\prime}}(\pi_r^t)=\{b^{\prime\prime}\in B^{**}:~~~\pi^{t***t}_r(b^{\prime\prime},a^{\prime\prime})=
\pi^{***}_r(b^{\prime\prime},a^{\prime\prime})\}.$$\\
It is clear that ~~~~~~~$$\bigcap_{a^{\prime\prime}\in A^{**}}Z_{a^{\prime\prime}}^t(B^{**})=Z_{A^{**}}^t(B^{**})=
Z(\pi_\ell^t),$$ ~~~~~~~ ~~\\ $$~~~~~~~~~~~~~~~~~~~~~~~~~~~~\bigcap_{a^{\prime\prime}\in A^{**}}Z_{a^{\prime\prime}}(B^{**})=Z_{A^{**}}(B^{**})=
Z(\pi_r).$$\\
The definition of $Z_{b^{\prime\prime}}^t(A^{**})$ and $Z_{b^{\prime\prime}}(A^{**})$ for some $b^{\prime\prime}\in B^{**}$ are the same.\\\\

\noindent{\it{\bf Theorem 2-13.}} Let $B$ be a Banach left $A-module$ and $A$ has a $LBAI$ $(e_{\alpha})_{\alpha}\subseteq A$ such that $e_{\alpha} \stackrel{w^*} {\rightarrow}e^{\prime\prime}$ in $A^{**}$ where $e^{\prime\prime}$ is a left unit for $A^{**}$. Suppose that $Z^t_{e^{\prime\prime}}(B^{**})=B^{**}$. {Then, $B$ factors on the right with respect to $A$ if and only if
$e^{\prime\prime}$ is a left unit for $B^{**}$.
\begin{proof}
Assume that  $B$ factors on the right with respect to $A$. Then  for every $b\in B$, there are $x\in B$ and $a\in A$ such that $b= ax$. Then for every $b^\prime\in B^*$, we have
$$<\pi^*_\ell(b^\prime,e_\alpha),b>=<b^\prime,\pi_\ell(e_\alpha,b)>=<\pi^{***}_\ell(e_\alpha,b),b^\prime>$$$$=
<\pi^{***}_\ell(e_\alpha,ax),b^\prime>=
<\pi^{***}_\ell(e_\alpha a,x),b^\prime>$$$$=<e_\alpha a,\pi^{**}_\ell(x,b^\prime)>=
<\pi^{**}_\ell(x,b^\prime),e_\alpha a>$$$$\rightarrow < \pi^{**}_\ell(x,b^\prime),a>=<b^\prime,b>.$$
It follows that  $\pi^*_\ell(b^\prime,e_\alpha) \stackrel{w^*} {\rightarrow}b^\prime$ in $B^*$.
Let $b^{\prime\prime}\in B^{**}$ and $(b_\beta)_\beta\subseteq B$ such that $b_\beta\stackrel{w^*} {\rightarrow}b^{\prime\prime}$ in $B^{**}$. Since $Z^t_{e^{\prime\prime}}(B^{**})=B^{**}$, for every $b^\prime\in B^*$, we have the following equality
$$<\pi_\ell^{***}(e^{\prime\prime},b^{\prime\prime}),b^\prime>=
\lim_\alpha\lim_\beta<b^\prime, \pi_\ell(e_\alpha,b_\beta)>$$$$=
\lim_\beta\lim_\alpha<b^\prime, \pi_\ell(e_\alpha,b_\beta)>=
\lim_\beta<b^\prime,b_\beta>$$$$
=<b^{\prime\prime},b^\prime>.$$
It follows that $\pi_\ell^{***}(e^{\prime\prime},b^{\prime\prime})=b^{\prime\prime}$, and so $e^{\prime\prime}$ is a left unit for $B^{**}$.\\
Conversely, let $e^{\prime\prime}$ be a left unit for $B^{**}$ and suppose that $b\in B$. Thren for every $b^\prime\in B^*$, we have
$$<b^\prime,\pi(e_\alpha,b)>=<\pi^{***}(e_\alpha,b),b^\prime>=<e_\alpha,\pi^{**}(b,b^\prime)>
=<\pi^{**}(b,b^\prime),e_\alpha>$$$$=
<e^{\prime\prime},\pi^{**}(b,b^\prime)>=<\pi^{***}(e^{\prime\prime},b),b^\prime>=<b^\prime,b>.$$
Then we have $\pi^*_\ell(b^\prime,e_\alpha) \stackrel{w} {\rightarrow}b^\prime$ in $B^*$, and so by Cohen factorization theorem we are done.

\end{proof}

\noindent{\it{\bf Corollary 2-14.}} Let $B$ be a Banach left $A-module$ and $A$ has a $LBAI$ $(e_{\alpha})_{\alpha}\subseteq A$ such that $e_{\alpha} \stackrel{w^*} {\rightarrow}e^{\prime\prime}$ in $A^{**}$ where $e^{\prime\prime}$ is a left unit for $A^{**}$. Suppose that $Z^t_{e^{\prime\prime}}(B^{**})=B^{**}$. Then
$\pi^*_\ell(b^\prime,e_\alpha) \stackrel{w} {\rightarrow}b^\prime$ in $B^*$ if and only if $e^{\prime\prime}$ is a left unit for $B^{**}$.\\\\

 For a Banach algebra $A$, we recall that a bounded linear operator $T:A\rightarrow A$ is said to be a left (resp. right) multiplier
if, for all $a, b\in A$, $T(ab)=T(a)b$ (resp. $T(ab)=aT(b)$). We denote by $LM(A)$ (resp. $RM(A)$) the set of all left (resp. right) multipliers of $A$. The set $LM(A)$ (resp. $RM(A)$) is normed subalgebra of the algebra $L(A)$ of bounded linear operator on $A$.\\
 Let  $B$ be a Banach left [resp. right] $A-module$ and  $T\in \mathbf{B}(A,B)$. Then $T$ is called extended left [resp. right] multiplier if ~$T(a_1 a_2)=\pi_r(T(a_1),a_2)\\ ~[resp.~T(a_1 a_2)=\pi_\ell(a_1,T(a_2))]$ ~for all $a_1,a_2\in A$.\\
We show by $LM(A,B)~[resp. ~ RM(A,B)]$ all of the Left [resp. right] multiplier extension from $A$ into $B$. \\
 Let $a^\prime\in A^*$. Then the mapping $T_{a^\prime}:a\rightarrow {a^\prime}a ~[resp.~~ R_{a^\prime}~a\rightarrow a{a^\prime}]$ from $A$ into $A^*$ is left [right] multiplier, that is, $T_{a^\prime}\in LM(A,A^*)~[R_{a^\prime}\in RM(A,A^*)].$ $T_{a^\prime}$ is weakly compact if and only if ${a^\prime}\in wap(A)$. So, we can write $wap(A)$ as a subspace of $LM(A,A^*)$.\\

\noindent{\it{\bf Theorem 2-15.}} Let $B$ be a   Banach $A-bimodule$ with a $BAI$ $(e_\alpha)_\alpha \subseteq A$. Then \begin{enumerate}
\item If $T\in LM(A,B)$, then $T(a)=\pi_r^{***}(b^{\prime\prime},a)$ for some $b^{\prime\prime}\in B^{**}$.
\item If $T\in RM(A,B)$, then $T(a)=\pi_\ell^{***}(a,b^{\prime\prime})$ for some $b^{\prime\prime}\in B^{**}$.
\end{enumerate}
\begin{proof}
\begin{enumerate}
 \item Since $(T(e_\alpha))_\alpha \subseteq B$ is bounded, it has weakly limit point in $B^{**}$. Let $b^{\prime\prime}\in B^{**}$ be a weakly limit point of $(T(e_\alpha))_\alpha$ and without loss generally, take  $T(e_\alpha) \stackrel{w} {\rightarrow}b^{\prime\prime}$. Then for every $b^\prime\in B^*$ and $a\in A$, we have

$$<\pi_r^{***}(b^{\prime\prime},a),b^\prime>
=\lim_\alpha<b^\prime, T(e_\alpha)a>=\lim_\alpha<b^\prime, T(e_\alpha a)>$$$$=\lim_\alpha<T^*(b^\prime), e_\alpha a>=
<T^*(b^\prime),  a>=<b^\prime,  T(a)>.$$
It follows that $\pi_r^{***}(b^{\prime\prime},a)=T(a)$.
\item  Proof is similar to (1).\\
\end{enumerate}\end{proof}

In the proceeding theorem, if we take $B=A$, then we have the following statements
\begin{enumerate}
\item If $T\in LM(A)$, then $T(a)=a^{\prime\prime}a$  for some $a^{\prime\prime}\in A^{**}$.
\item If $T\in RM(A)$, then $T(a)=aa^{\prime\prime}$  for some $a^{\prime\prime}\in A^{**}$.\\
\end{enumerate}

 Let  $B$ be a  Banach left [resp. right] $A-module$. Then for every $b\in B$, we define $\ell_b$ (resp. $r_b$) the linear mapping $a\rightarrow \pi_\ell(a,b)$ (resp. $a\rightarrow \pi_r(b,a)$).\\
In the proceeding theorem, for $T\in LM(A,B)~[resp. ~ T\in RM(A,B)]$ if $B$ is weakly compact, then   $T=r_b$ [resp. $T=\ell_b$], for some $b\in B$.\\\\
Let $B$ be a   Banach right $A-module$. Then, we define $B^*B$ as a subspace of $A^*$ including of all $\pi_r(b^\prime ,b)$ for every $b^\prime \in B^*$ and $b\in B$, that is, for every $a\in A$, we define
$$<\pi^*_r(b^\prime ,b), a>=<b^\prime ,\pi_r(b, a)>.$$
 For Banach left $A-module$ $B$, we also define $B^{**}B^*$ as a subset of $A^*$ and for every $b^{\prime\prime}\in B^{**}$, $b^\prime \in B^*$ and $a\in A$, as follows
$$<\pi^{**}_\ell(b^{\prime\prime},b^\prime ),a>=<b^{\prime\prime},\pi^*_\ell(b^\prime ,a)>.$$\\
A Banach space $B$ is said to be weakly
sequentially complete ($WSC$), if every weakly Cauchy sequence in $B$ has
a weak limit in $B$.\\\\

\noindent{\it{\bf Definition 2-16.}} Let $B$ be a Banach left  $A-module$ and $b^{\prime\prime}\in B^{**}$. Suppose that $(b_{\alpha})_{\alpha}\subseteq B$ such that  $b_{\alpha} \stackrel{w^*} {\rightarrow}b^{\prime\prime}$. We define the following set
$$\widetilde{Z}_{b^{\prime\prime}}(A^{**})=\{a^{\prime\prime}\in A^{**}:~~~\pi_\ell^{***}(a^{\prime\prime},b_{\alpha})\stackrel{w^*} {\rightarrow}\pi_\ell^{***}(a^{\prime\prime},b^{\prime\prime})\},$$
which is subspace of $A^{**}$. It is clear that $Z_{b^{\prime\prime}}(A^{**})\subseteq \widetilde{Z}_{b^{\prime\prime}}(A^{**})$, and so
$$Z_{B^{**}}(A^{**})=\bigcap_{b^{\prime\prime}\in B^{**}}{Z}_{b^{\prime\prime}}(A^{**})\subseteq \bigcap_{b^{\prime\prime}\in B^{**}}\widetilde{Z}_{b^{\prime\prime}}(A^{**}).$$
For a Banach right $A-module$, the definition of $\widetilde{Z}^t_{a^{\prime\prime}}(B^{**})$  is similar.\\\\

Let $B$ be a Banach left $A-module$. Then we define $B^{**}B^*$ as follows.
$$B^{**}B^*=\{\pi_\ell^{**}(b^{\prime\prime},b^{\prime}):~b^{\prime\prime}\in B^{**} ~ and~b^\prime\in B^*\}.$$
 It is clear that $B^{**}B^*$ is a subspace of $A^*$.\\\\

Let $A$ be a Banach algebra. Then,  $A$ is said to be weakly
sequentially complete ($=WSC$), if every weakly Cauchy sequence in $A$ has
a weak limit.\\\\

\noindent{\it{\bf Theorem 2-17.}} Let $B$ be a left  Banach $A-module$ and $T\in \mathbf{B}(A,B)$. Consider the following statements.
\begin{enumerate}
\item $T=\ell_b$, for some $b\in B$.
\item  $T^{**}(a^{\prime\prime})=\pi_\ell^{***}(a^{\prime\prime},b^{\prime\prime})$ for some $b^{\prime\prime}\in B^{**}$ such that $\widetilde{Z}_{b^{\prime\prime}}(A^{**})=A^{**}$.
\item $T^*(B^*)\subseteq BB^*$.
\end{enumerate}
Then $(1)\Rightarrow (2)\Rightarrow (3)$.\\
Assume that $B$ has $WSC$. If we take $T\in RM(A,B)$ and  $B$ has a sequential $BAI$, then (1), (2) and (3) are equivalent.\\
\begin{proof}
$(1)\Rightarrow (2)$\\
 Let $T=\ell_b$, for some $b\in B$. Then $T^{**}(a^{\prime\prime})=\ell^{**}_b(a^{\prime\prime})=\pi_\ell^{***}(a^{\prime\prime},b)$ for every $a^{\prime\prime}\in A^{**}$, and so proof is hold.\\
 $(2)\Leftrightarrow (3)$\\
Take $a^{\prime\prime}\in (BB^*)^\bot$. Assume that $b^{\prime\prime}\in B^{**}$ and $(b_{\alpha})_{\alpha}\subseteq B$ such that  $b_{\alpha} \stackrel{w^*} {\rightarrow}b^{\prime\prime}$. For every $b^\prime\in B^{**}$, we have the following equality
$$<a^{\prime\prime}, T^*(b^\prime)>=<T^{**}(a^{\prime\prime}), b^\prime>=<\pi_\ell^{***}(a^{\prime\prime},b^{\prime\prime}), b^\prime>=\lim_\alpha<\pi_\ell^{***}(a^{\prime\prime},b_{\alpha}),b^{\prime}>$$$$=
\lim_\alpha<a^{\prime\prime},\pi_\ell^{**}(b_{\alpha},b^{\prime})>=0.$$
It follows that $T^*(B^*)\subseteq BB^*$.\\
Take $T\in RM(A,B)$ and suppose that  $B$ is $WSC$ with sequential $BAI$. It is enough, we show that $(3)\Rightarrow (1)$.
Assume that $(e_n)_n\subseteq A$ is a $BAI$ for $B$. Then for every $b^\prime\in B^*$, we have
$$\mid<b^\prime,T(e_n)>-<b^\prime,T(e_m)>\mid=\mid <T^*(b^\prime),e_n-e_m>\mid=
\mid <\pi^{**}_\ell(b,b^\prime),e_n-e_m>\mid$$$$=
\mid <b,\pi^{*}_\ell(b^\prime,e_n-e_m)>\mid=
\mid <b^\prime,\pi^{}_\ell(e_n-e_m,b)>\mid\rightarrow 0.$$
It follows that $(T(e_n))_n$ is weakly Cauchy sequence in $B$ and since $B$ is $WSC$, there is $b\in B$ such that
$T(e_n)\stackrel{w} {\rightarrow}b$ in $B$. Let $a\in A$. Then for every $b^\prime\in B^*$, we have
$$<b^\prime ,\pi_\ell(a,b)>=<\pi_\ell^*(b^\prime ,a),b)>=\lim_n<\pi_\ell^*(b^\prime ,a),T(e_n)>$$$$
=\lim_n<b^\prime ,\pi_\ell(a,T(e_n))>=\lim_n<b^\prime ,T(ae_n)>$$$$=\lim_n<T^*(b^\prime),ae_n>
=<T^*(b^\prime),a>$$$$=<b^\prime,T(a)>.$$
Thus $\ell_b(a)=\pi_\ell(a,b)=T(a)$.\\
\end{proof}

\noindent{\it{\bf Example 2-18.}} Let $G$ be a locally compact group. Then by convolution multiplication,  $M(G)$ is a  $L^1(G)-bimodule$. Let $f\in L^1(G)$ and $T(\mu)=\mu*f$ for all $\mu\in M(G)$. Then $T^*(L^\infty (G))\subseteq M(G)M(G)^*$. Also if we take  $T(\mu)=f*\mu$ for all $\mu\in M(G)$, then we have $T^*(L^\infty (G))\subseteq M(G)^*M(G)$.

\bibliographystyle{amsplain}

\it{Department of Mathematics, Amirkabir University of Technology, Tehran, Iran\\
{\it Email address:} haghnejad@aut.ac.ir\\

\end{document}